\documentclass[fontsize=12pt]{scrartcl}

\usepackage[british]{babel}
\usepackage{lmodern}
\usepackage[T1]{fontenc}
\usepackage{textcomp}
\usepackage{amsmath}
\usepackage{amsfonts}
\usepackage{amssymb}
\usepackage{enumerate}
\usepackage[amsmath, thmmarks]{ntheorem}
\usepackage{paralist}
\usepackage{booktabs}

\addtokomafont{disposition}{\rmfamily}

\theorembodyfont{\rmfamily}
\theoremstyle{change}
\theoremseparator{.}
\theoremsymbol{\ensuremath{\diamond}}
\newtheorem{thmetc}{Dummy}[section]
\newtheorem{theorem}[thmetc]{Theorem}
\newtheorem{lemma}[thmetc]{Lemma}
\newtheorem{corollary}[thmetc]{Corollary}
\newtheorem{example}[thmetc]{Example}
\newtheorem{examples}[thmetc]{Examples}
\newtheorem{remark}[thmetc]{Remark}
\newtheorem{algorithm}[thmetc]{Algorithm}
\theoremstyle{nonumberplain}
\theoremsymbol{\ensuremath{\square}}
\newtheorem{proof}{Proof}

\newcommand{\N}{\mathbb{N}}
\newcommand{\Z}{\mathbb{Z}}

\newcommand{\GL}{\operatorname{GL}}
\newcommand{\NN}{\operatorname{N}}
\renewcommand{\a}{\mathfrak{a}}
\renewcommand{\b}{\mathfrak{b}}
\renewcommand{\c}{\mathfrak{c}}
\renewcommand{\d}{\mathfrak{d}}
\newcommand{\e}{\mathfrak{e}}
\newcommand{\f}{\mathfrak{f}}
\newcommand{\g}{\mathfrak{g}}
\renewcommand{\o}{\mathfrak{o}}
\newcommand{\p}{\mathfrak{p}}
\newcommand{\q}{\mathfrak{q}}
\renewcommand{\r}{\mathfrak{r}}
\renewcommand{\v}{\operatorname{v}}
\newcommand{\MM}{\mathfrak{M}}
\newcommand{\QQ}{\mathfrak{Q}}
\newcommand{\abs}[1]{\left|#1\right|}
\newcommand{\tabs}[1]{\mathopen|#1\mathclose|}

\newcommand{\matr}[4]{\begin{pmatrix} #1 & #2 \\ #3 & #4 \end{pmatrix}}
\newcommand{\tmatr}[4]{(\begin{smallmatrix} #1 & #2 \\ #3 & #4 \end{smallmatrix})}

\newcommand{\tvect}[2]{(\begin{smallmatrix} #1 \\ #2 \end{smallmatrix})}
\newcommand{\divides}{\mathrel|}
\newcommand{\ndivides}{\mathrel\nmid}
\newcommand{\eps}{\varepsilon}

\title{Counting cosets of unimodular groups over Dedekind domains}
\author{Marc Ensenbach\textsuperscript{a}\footnote{\texttt{ensenbach@mathematik.uni-siegen.de}, Telephone: +49\,271\,740\,3518, Fax: +49\,271\,740\,3514} \\ \normalsize \textsuperscript{a}University of Siegen, Department of Mathematics, 57068 Siegen, Germany}
\date{}
\begin{document}
\maketitle

\begin{abstract}
In this paper, a formula for the calculation of the number of right cosets contained in a double coset with respect to a unimodular group over a Dedekind domain is developed, and applications of this formula in the theory of congruence subgroups -- an index formula -- and the theory of abstract Hecke algebras -- a reduction theorem and an algorithm for the explicit calculation of products -- are given.
\end{abstract}

\emph{Keywords:} Unimodular group, Dedekind domain, Congruence subgroup, Index formula, Hecke algebra

\emph{2010 MSC:} 20G30, 20H05, 20C08

\section{Introduction}\label{intro}

Being able to count right cosets contained in a double coset with respect to a unimodular group $U$, i.\,e.\ being able to determine the cardinality of $U \backslash U A U$ for any given matrix $A$ of appropriate format, is helpful in a multitude of areas.

The formula presented in this article originally has been studied in the field of Hecke algebras. If we want to carry out computational analysis in Hecke theory, we need an algorithm that allows us to multiply two elements of an abstract Hecke algebra. Since the product can be calculated by a multiplication of representatives of right cosets, the task of multiplying elements of an abstract Hecke algebra can essentially be reduced to the search for decompositions of double cosets into right cosets. Knowing how many cosets we have to find allows us to state a randomised algorithm which carries out the decomposition. Moreover, there is an application in the proof of a theoretical result on abstract Hecke algebras. In the ``classic'' Hecke algebra $H_n$ related to the $\GL_n(\Z)$, certain products in $H_n$ can be reduced to products in $H_{n - 1}$. With the presented result this reduction theorem can be generalised to Hecke algebras related to arbitrary norm-finite Dedekind domains in the case $n = 2$. Details of these two applications to Hecke algebras can be found in section \ref{applhecke}.

Another application is shown in section \ref{applcong}. The set of right cosets with respect to $U$ which are contained in $U A U$ can be related to a right transversal of a certain subgroup of $U$ depending on $A$. For some special matrices $A$ these subgroups turn out to be congruence subgroups, so we are enabled to calculate indexes of certain congruence subgroups of $U$.

The remaining sections of this article are organised as follows: In the following section, the notation used in this article is fixed and some basic facts which are used throughout this article are assembled. After that, the main result -- a formula for the cardinality of $U \backslash U A U$ -- is stated and proved in section \ref{countright}.

This article is developed from a talk I gave some time ago in the research seminar ``Computational Algebra and Number Theory'' of Fritz Grunewald, whose unexpected death in 2010 means a great loss for the mathematical community.

\section{Preliminaries and Notation}\label{prelim}

Denote by $\o$ a norm-finite Dedekind domain, i.\,e.\ a Dedekind domain in which $\tabs{\o / a \o} < \infty$ holds for every $a \in \o$ (where $\tabs{M}$ is the cardinality of the set $M$). Furthermore, denote by $K$ the field of fractions of $\o$, and by $\o^*$ the group of unities of $\o$; then denote by $\v_\p(\a)$ the multiplicity of a prime ideal $\p$ in the ideal $\a$ of $\o$ (fundamental properties of Dedekind domains and multiplicities can be found for example in \cite{Froehlich:Algebraic} Chapter II).

Let $I$ be the set of ($2 \times 2$) matrices with entries in $\o$ and non-zero determinant; furthermore, denote by $U$ the set of matrices in $I$ with determinant in $\o^*$ (in other words $U = \GL_2(\o)$ and $I = \GL_2(K) \cap \o^{2 \times 2}$). For $A = \tmatr{a}{b}{c}{d} \in \o^{2 \times 2}$ one defines the first and second determinantal divisor of $A$ by $\d_1(A) = a \o + b \o + c \o + d \o$ and $\d_2(A) = (\det A) \o$, respectively. Furthermore, $\e_1(A) := \d_1(A)$ and $\e_2(A) := \d_2(A) \d_1(A)^{-1}$ are called the elementary divisors of $A$. Additionally, define the fundamental factors $\f_1(A) := \e_1(A) = \d_1(A)$ and $\f_2(A) := \e_2(A) \e_1(A)^{-1} = \d_2(A) \d_1(A)^{-2}$ and introduce the notation $\g(A)$ for the g.\,c.\,d.\ of the first column of $A$, i.\,e.\ $\g(A) = a \o + c \o$.

The relation between determinantal divisors and double cosets of $U$ is given in the following theorem, which goes back to Steinitz (\cite{Steinitz:Rechteckige}, see also \cite{Ensenbach:Determinantal} Theorem 2.2).

\begin{theorem}\label{dcequal}
Let $A, B \in \o^{2 \times 2}$.
\begin{enumerate}[a)]
\item If $A$ and $B$ have rank $2$ (i.\,e., if $A, B \in I$), the following assertions are equivalent:
\begin{inparaenum}[(i)]
\item $U A U = U B U$,
\item $\d_1(A) = \d_1(B)$ and $\d_2(A) = \d_2(B)$.
\end{inparaenum}
\item If $A$ and $B$ have rank $1$ and the first columns of $A$ and $B$ both contain at least one non-zero element, the following assertions are equivalent:
\begin{inparaenum}[(i)]
\item $U A U = U B U$,
\item $\d_1(A) = \d_1(B)$ and $\g(A) = \g(B)$.
\end{inparaenum}
\end{enumerate}
In these assertions, the $\d_i$ can also be replaced by the $\e_i$ or the $\f_i$.
\end{theorem}

This theorem can not only be used to characterise the equality of double cosets, but also has an application in the proof of the following corollary which allows to state a relation between different generators of the same ideal in $\o$.

\begin{corollary}\label{unimodgen}
Let $a, b, c, d \in \o$ such that $a \o + b \o = c \o + d \o$. Then there exists an $R \in U$ satisfying $R \tvect{a}{b} = \tvect{c}{d}$.
\end{corollary}
\begin{proof}
In the case $c = d = 0$ we also have $a = b = 0$ and can choose $R = \tmatr{1}{0}{0}{1}$. For the remaining part of the proof assume $c \neq 0$ (without loss of generality). Let $A = \tmatr{a}{0}{b}{0}$ and $B = \tmatr{c}{0}{d}{0}$. Since $A$ and $B$ both have rank $1$ and satisfy $\d_1(A) = \d_1(B)$ as well as $\g(A) = \g(B)$, Theorem \ref{dcequal} yields the existence of $P, Q \in U$ such that $P A Q = B$ and thus $P A = B Q^{-1}$. Writing $P = \tmatr{p_1}{p_2}{p_3}{p_4}$ and $Q^{-1} = \tmatr{q_1}{q_2}{q_3}{q_4}$ and calculating $P A$ as well as $B Q^{-1}$ we obtain
\[
\matr{p_1 a + p_2 b}{0}{p_3 a + p_4 b}{0} = \matr{c q_1}{c q_2}{d q_1}{d q_2}.
\]
In particular, we have $c q_2 = 0$, and since $c \neq 0$, this implies $q_2 = 0$. Thus $\det Q^{-1} = q_1 q_4$, which implies $q_1 \in \o^*$. If we define $R = q_1^{-1} P$, we thus have $R \in U$. Furthermore,
\[
R \matr{a}{0}{b}{0} = q_1^{-1} P A = q_1^{-1} \matr{p_1 a + p_2 b}{0}{p_3 a + p_4 b}{0} = q_1^{-1} \matr{c q_1}{0}{d q_1}{0} = \matr{c}{0}{d}{0},
\]
which proves $R \tvect{a}{b} = \tvect{c}{d}$ and completes the proof.
\end{proof}

\begin{remark}
Since there exists a version of Theorem \ref{dcequal} for $A, B \in \o^{n \times n}$ for arbitrary $n \in \N$ (see e.\,g.\ \cite{Ensenbach:Determinantal} Theorem 2.2), Corollary \ref{unimodgen} can easily be generalised from two generators to an arbitrary number of generators of an ideal, as long as the number of generators on both sides of the equation are the same.
\end{remark}

\section{Counting right cosets}\label{countright}

In this section a formula for the number of right cosets in a given double coset of $U$ is derived. To begin with, a short example shows how coset counting is carried out in the ``nice classic case'' $\o = \Z$. This will serve as a guideline for the subsequent analysis of the general case.

\begin{example}\label{exz}
Let $\o = \Z$. Since $\o$ is a principal ideal domain, every right coset $U B$ for $B \in I$ has a unique representative
\[
B' = \matr{a}{b}{0}{d} \qquad\text{with $a, d > 0$ and $0 \leq b < d$}
\]
known as the Hermite normal form of $B$. Given this normal form, the number of right cosets in a given double coset $U A U$ can be obtained by generating all possible normal forms (in a sensible way) and deciding whether they belong to $U A U$. The latter can be carried out using Theorem \ref{dcequal} to test for $U B' U = U A U$, so it has to be checked whether $\d_1(B') = \d_1(A)$ and $\d_2(B') = \d_2(A)$ hold.

As a concrete example construct every right coset representative $B'$ as above contained in $U A U$ where $A = \tmatr{1}{0}{0}{4}$. Since $\d_2(A) = \d_2(B')$ is a necessary condition for $U B' U = U A U$, the equation $(\det A) \Z = (\det B') \Z$ and thus $4 = a d$ has to be satisfied. So there are three possible cases:
\begin{inparaenum}[(i)]
\item $a = 4$ and $d = 1$,
\item $a = 2$ and $d = 2$, and
\item $a = 1$ and $d = 4$.
\end{inparaenum}
For these cases determine, for which values of $b$ the equation $U B' U = U A U$ is fulfilled. To this end, it suffices to test whether $\d_1(B') = \d_1(A)$ holds since $a$ and $d$ have already been constructed to satisfy $\d_2(B') = \d_2(A)$.
\begin{asparaenum}[{Case} (i):]
\item Since $d = 1$ and $0 \leq b < d = 1$, only the case $b = 0$ has to be analysed. Then we have $\d_1(B') = a \Z + b \Z + d \Z = 4 \Z + 0 \Z + 1 \Z = 1 \Z = 1 \Z + 4 \Z = \d_1(A)$, so $\tmatr{4}{0}{0}{1}$ is a right coset representative in $U A U$.
\item Since $d = 2$ and $0 \leq b < d = 2$, the cases $b = 0$ and $b = 1$ have to be considered. For $b = 0$ we have $\d_1(B') = 2 \Z \neq 1 \Z = \d_1(A)$, so $\tmatr{2}{0}{0}{2}$ is not an element of $U A U$. For $b = 1$, however, we have $\d_1(B') = 1 \Z = \d_1(A)$, so $\tmatr{2}{1}{0}{2}$ belongs to $U A U$.
\item Since $d = 4$ and $0 \leq b < d = 4$, the cases $b \in \{0, 1, 2, 3\}$ have to be analysed. Due to $a = 1$ we have $\d_1(B') = 1 \Z = \d_1(A)$ in any of these cases, so $\tmatr{1}{b}{0}{4}$ belongs to $U A U$ for every $b \in \{0, 1, 2, 3\}$.
\end{asparaenum}

Summarising, in $U A U$ we have found the $6$ right coset representatives $\tmatr{4}{0}{0}{1}$, $\tmatr{2}{1}{0}{2}$, $\tmatr{1}{0}{0}{4}$, $\tmatr{1}{1}{0}{4}$, $\tmatr{1}{2}{0}{4}$, and $\tmatr{1}{3}{0}{4}$.

Generalising these considerations, a formula for the number $\mu(A)$ of right cosets contained in $U A U$ can be stated:
\[
\mu(A) = \sum_{\substack{d \in \N \\ d \divides \det A}} \tabs{\{b \in \N_0 \:|\: b < d \text{ and } \tfrac{\det A}{d} \Z + b \Z + d \Z = \d_1(A)\}}
\]
(where $\N = \{1, 2, 3, \ldots, \}$ and $\N_0 = \{0, 1, 2, \ldots\}$). The cardinality of a set $\{b \in \N_0 \:|\: b < d \text{ and } a \Z + b \Z + d \Z = \d_1(A)\}$ can be calculated explicitly, which finally leads to a product formula for $\mu(A)$ (not presented in detail since the same steps are to be done for the general case in the following).
\end{example}

The first main ingredient of the approach taken in Example \ref{exz} in the classic case was the Hermite normal form. In the general case, another normal form can be constructed -- not as ``nice'' as in the classic case, but nevertheless solving the issue of a uniquely determined representative.

\begin{lemma}\label{normalform}
Let $\a$ be an ideal in $\o$ and $b \in \a$. Choose an $a \in \a$ satisfying $a \neq 0$ and $a \o + b \o = \a$ (always possible since $\o$ is a Dedekind domain) and a transversal $T$ of $(\o \cap a b^{-1} \o) / (b a^{-1} \o \cap a b^{-1} \o)$. This transversal is finite, and for every $A \in I$ satisfying $\g(A) = \a$ and $\d_2(A) = b \o$ there exists a uniquely determined $c \in T$ such that
\[
U \matr{a}{c - 1}{b}{b a^{-1} c} = U A.
\]
\end{lemma}
\begin{proof}
The finiteness of $T$ follows from the norm-finiteness of $\o$ since $b a^{-1} \o \cap a b^{-1} \o$ is an ideal in $\o$.

To prove the uniqueness of $c$, assume
\[
U \matr{a}{c - 1}{b}{b a^{-1} c} = U \matr{a}{d - 1}{b}{b a^{-1} d}
\]
for $c, d \in T$. Since then
\[
\matr{a}{d - 1}{b}{b a^{-1} d} \matr{a}{c - 1}{b}{b a^{-1} c}^{-1} = \matr{c - d + 1}{a b^{-1} (d - c)}{b a^{-1} (c - d)}{1 - c + d}
\]
has to be an element of $U$, we obtain in particular $a b^{-1} (d - c) \in \o$ and $b a^{-1} (c - d) \in \o$, which yields $d - c \in b a^{-1} \o \cap a b^{-1} \o$. Since $T$ is a transversal modulo $b a^{-1} \o \cap a b^{-1} \o$ and $c, d \in T$, this shows $c = d$ and thus proves the uniqueness of the representative.

In the remaining part of the proof the existence of the desired representative is shown. Let $A \in I$ satisfying $\g(A) = \a$ and $\d_2(A) = b \o$. Since $a \o + b \o = \g(A)$, by Corollary \ref{unimodgen} there exists a $P_1 \in U$ such that $P_1 A = \tmatr{a}{*}{b}{*}$. Then let $\eps = b (\det(P_1 A))^{-1} \in \o^*$ and $P_2 = \tmatr{\eps}{0}{0}{1} \in U$, such that $\det(P_2 P_1 A) = b$. Furthermore, the Chinese Remainder Theorem allows us to choose a $p \in \o$ satisfying $p \in (\o \cap b a^{-1} \o) + \eps^{-1}$ and $p \in (\o \cap a b^{-1} \o) + 1$ since $\o \cap b a^{-1}$ and $\o \cap a b^{-1} \o$ are relatively prime. The matrix
\[
P_3 = \matr{p}{a b^{-1} (1 - \eps p)}{b a^{-1} (p - 1)}{\eps + 1 - \eps p}
\]
then is an element of $\o^{2 \times 2}$ with $\det P_3 = 1$, and we have
\[
P_3 P_2 P_1 A = P_3 P_2 \matr{a}{*}{b}{*} = P_3 \matr{\eps a}{*}{b}{*} = \matr{a}{*}{b}{*}
\]
with $\det (P_3 P_2 P_1 A) = b$, so if the second column of $P_3 P_2 P_1 A$ is denoted by $\tvect{r}{s}$, we have $a s - b r = b$ and thus $1 + r = a b^{-1} s \in a b^{-1} \o$. Since furthermore $1 + r \in \o$, by the choice of $T$ there exists a $c \in T$ satisfying $1 + r \in c + (b a^{-1} \o \cap a b^{-1} \o)$. Now let
\[
P_4 = \matr{a b^{-1} s - c + 1}{a b^{-1} (c - r - 1)}{s - b a^{-1} c}{c - r}.
\]
Then $P_4 \in \o^{2 \times 2}$ and $\det P_4 = 1$ (since $a b^{-1} s = 1 + r$ and $c - r - 1 \in b a^{-1} \o \cap a b^{-1} \o$), and putting everything together we have $P_4 P_3 P_2 P_1 \in U$ and
\[
P_4 P_3 P_2 P_1 A = P_4 \matr{a}{r}{b}{s} = 
\matr{a (a b^{-1} s - r)}{a b^{-1} s (c - 1) - c r + r}{a s - b r}{b a^{-1} c (a b^{-1} s - r)}
=
\matr{a}{c - 1}{b}{b a^{-1} c},
\]
which shows the existence of a representative with the desired form and thus completes the proof.
\end{proof}

With this normal form, a first elementary formula for the number of right cosets with a prescribed g.\,c.\,d.\ of the first column can be given. (All elements of a right coset with respect to $U$ have the same g.\,c.\,d.\ of the first column, so it is possible to talk about the g.\,c.\,d.\ of the first column of a right coset.)

\begin{corollary}\label{elemmua}
Let $A \in I$ and $\a$ an ideal of $\o$ such that $\d_1(A) \divides \a \divides \d_2(A)$. Choose an $a \in \a$ satisfying $a \o + \d_2(A) = \a$ (possible since $\a \divides \d_2(A)$ and $\o$ is a Dedekind domain), let $\q = \a^{-1} a$ as well as $\b = \a^{-1} \d_2(A)$, and choose a transversal $T$ of $\q / \q \b$. Then the number $\mu_\a(A)$ of right cosets in $U A U$ with $\a$ as g.\,c.\,d.\ of the first column can be calculated by
\[
\mu_\a(A) = \tabs{\{c \in T \:|\: \a + (c - 1) \o + c \q^{-1} \b = \d_1(A)\}}.
\]
\end{corollary}
\begin{proof}
Since $\a = a \o + \d_2(A) = \a \q + \a \b$, the ideals $\q$ and $\b$ are relatively prime, which yields $\o \cap a \d_2(A)^{-1} = \o \cap \q \b^{-1} = \q$ as well as $\d_2(A) a^{-1} \cap a \d_2(A)^{-1} = \b \q^{-1} \cap \q \b^{-1} = \b \q$. Thus $T$ is a transversal of $(\o \cap a \d_2(A)^{-1}) / (\d_2(A) a^{-1} \cap a \d_2(A)^{-1})$.

If $U B$ for some $B \in I$ is a right coset in $U A U$ satisfying $\g(B) = \a$, then in particular $\d_2(B) = \d_2(A)$, and according to Lemma \ref{normalform} there exists a uniquely determined representative $C$ of $U B$ of the form described in that Lemma (with $b = \det A$ and $c \in T$, the latter according to the first paragraph of this proof). Thus
\begin{align*}
& \{U B \:|\: B \in I \text{ with } \g(B) = \a \text{ and } U B \subseteq U A U\} \\
& = \left\{U B \:\middle|\: B = \matr{a}{c - 1}{\det A}{(\det A) a^{-1} c} \text{ for some } c \in T \text{ and } B \in U A U\right\},
\end{align*}
and since $\d_2(B) = \d_2(A)$ for those $B$, Theorem \ref{dcequal} and $a \o + \d_2(A) = \a$ yield
\begin{align*}
\mu_\a(A) & = \abs{\left\{c \in T \:\middle|\: \d_1\left(\matr{a}{c - 1}{\det A}{(\det A) a^{-1} c} \right) = \d_1(A)\right\}} \\
& = \tabs{\{c \in T \:|\: \a + (c - 1) \o + c \q^{-1} \b = \d_1(A)\}}.
\end{align*}

\end{proof}

The formula presented in Corollary \ref{elemmua} is only a first step since it is not very far from a mere enumeration of right cosets. The next step is the establishment of a product formula for the cardinality on the right-hand side. To achieve this, we first need some auxiliary results.

\begin{lemma}\label{chardiv}
In the setting of Corollary \ref{elemmua} the following assertions are equivalent:
\begin{enumerate}[(i)]
\item $\c \divides \a + (c - 1) \o + c \q^{-1} \b$.
\item $\c \divides \a$ and $\a \c \divides \d_2(A) $ and $\c \divides c - 1$.
\end{enumerate}
\end{lemma}
\begin{proof}
First assume that (i) is satisfied and show that (ii) is fulfilled. Since (i) implies $\c \divides (c - 1) \o + c \q^{-1} \b$ and $c$ and $c - 1$ are relatively prime, we have $\c \divides \q^{-1} \b$ and thus $\a \c \divides \a \q^{-1} \b = a^{-1} \a \d_2(A)$, which implies $\a \c \divides \d_2(A)$ since $a \in \a$. The remaining parts of (ii) follow obviously.

Now assume that (ii) is fulfilled. For the proof of (i) it remains to show $\c \divides c \q^{-1} \b$. But this follows from $\a \c \divides \d_2(A)$ and $c \in \q$ since the latter implies $\c \divides \a^{-1} \d_2(A) = \b$ and $c \q^{-1} \subseteq \o$, so the proof is complete.
\end{proof}

The second auxiliary result gives an explicit formula for the cardinality of a certain subset of $T$ needed in the calculation of $\mu_\a(A)$.

\begin{lemma}\label{tnorm}
In the setting of Corollary \ref{elemmua} we have
\[
\tabs{\{c \in T \:|\: c - 1 \in \c\}} = \NN(\b) \NN(\c)^{-1}
\]
for every ideal $\c$ of $\o$ satisfying $\c \divides \a \divides \d_2(A) \c^{-1}$.
\end{lemma}
\begin{proof}
In the given setting we have $\c \divides \b$, and $\b$ and $\q$ are relatively prime, so $\c$ and $\q$ are relatively prime. By the Chinese Remainder Theorem there exists a $d \in \q \cap (\c + 1)$. Then $\{c - d \:|\: c \in T,\; c - 1 \in \c\}$ is a transversal of $\c \q / \b \q$: All those $c - d$ are elements of $\c \cap \q = \c \q$, the $(c - d) + \b \q$ are pairwise different since $T$ is a transversal of $\q / \b \q$, and for every $x \in \c \q$ there exists a $c \in T$ satisfying $x \in (c - d) + \b \q$, namely the one satisfying $x + d \in c + \b \q$, which exists in $T$ since $x + d \in \q$ and is an element of $\{c \in T \:|\: c - 1 \in \c\}$ since $x + d \in c + \b \q \subseteq c + \c$ and thus $c - 1 \in x + (d - 1) + \c = \c$. Now the definition and the multiplicity of the norm yields $\tabs{\{c \in T \:|\: c - 1 \in \c\}} = \tabs{\c \q / \b \q} = \NN(\b) \NN(\c)^{-1}$.
\end{proof}

Now we are prepared to prove a product formula for $\mu_\a(A)$.

\begin{theorem}\label{tmua}
Let $A \in I$ and $\a$ be an ideal in $\o$. If $\d_1(A) \divides \a \divides \d_2(A) \d_1(A)^{-1}$, then
\[
\mu_\a(A) = \frac{\NN(\d_2(A))}{\NN(\a) \NN(\d_1(A))} \prod_{\substack{\p \text{ prime ideal} \\ \p \divides \a \d_1(A)^{-1} + \d_2(A) \d_1(A)^{-1} \a^{-1}}} (1 - \NN(\p)^{-1}),
\]
otherwise, $\mu_\a(A) = 0$ holds.
\end{theorem}
\begin{proof}
If $\mu_\a(A) > 0$, then there exists a $B \in U A U$ having $\a$ as g.\,c.\,d.\ of the first column. Then $\d_1(A) = \d_1(B) \divides \a$ and $\a \divides \d_2(B) = \d_2(A)$, so Corollary \ref{elemmua} is applicable. Thus, using the notation introduced in Corollary \ref{elemmua}, there exists a $c \in T$ satisfying $\a + (c - 1) \o + c \q^{-1} \b = \d_1(A)$. So Lemma \ref{chardiv} implies $\a \d_1(A) \divides \d_2(A)$, which shows that $\a \divides \d_2(A) \d_1(A)^{-1}$ is necessary for $\mu_\a(A) > 0$. Thus it is proved that $\mu_\a(A) = 0$ if $\d_1(A) \divides \a \divides \d_2(A) \d_1(A)^{-1}$ does not hold.

In the following assume $\d_1(A) \divides \a \divides \d_2(A) \d_1(A)^{-1}$. Denote by $\QQ$ the set of all prime ideals of $\o$ dividing $\d_2(A) \d_1(A)^{-2}$, let $M = \{c \in T \:|\: \a + (c - 1) \o + c \q^{-1} \b = \d_1(A)\}$ and $M(\c) = \{c \in T \:|\: \c \text{ divides } \a + (c - 1) \o + c \q^{-1} \b\}$ for all ideals $\c$ of $\o$. By the inclusion-exclusion principle we then have
\[
\tabs{M} = \sum_{\MM \subseteq \QQ} (-1)^{\tabs{\MM}} \abs{M\left(\d_1(A) \prod_{\q \in \MM} \q\right)}.
\]
If for a product $\q'$ of pairwise distinct prime ideals $\d_1(A) \q' \divides \a \divides \d_2(A) (\d_1(A) \q')^{-1}$ does not hold, then for $\c = \d_1(A) \q'$ the first or the second condition in Lemma \ref{chardiv}\,(ii) is violated, which implies $\tabs{M(\c)} = 0$. Since $\d_1(A) \q' \divides \a \divides \d_2(A) (\d_1(A) \q')^{-1}$ is equivalent to $\q' \divides \a \d_1(A)^{-1}$ and $\q' \divides \d_2(A) \d_1(A)^{-1} \a^{-1}$, in the above formula $\QQ$ can be replaced by $\QQ'$ where $\QQ'$ denotes the set of all prime ideals of $\o$ dividing $\a \d_1(A)^{-1} + \d_2(A) \d_1(A)^{-1} \a^{-1}$. In the case that $\q'$ is a product of pairwise distinct prime ideals in $\QQ'$, the condition $\d_1(A) \q' \divides \a \divides \d_2(A) (\d_1(A) \q')^{-1}$ is satisfied, and Lemma \ref{chardiv} and Lemma \ref{tnorm} yield $M(\d_1(A) \q') = \tabs{\{c \in T \:|\: \d_1(A) \q' \text{ divides } c - 1\}} = \NN(\b) \NN(\d_1(A) \q')^{-1}$. Plugging this into the above formula and using the multiplicity of the norm and the distributive law we obtain
\begin{align*}
\tabs{M} & = \frac{\NN(\b)}{\NN(\d_1(A))} \sum_{\MM \subseteq \QQ'} (-1)^{\tabs{\MM}} \NN\left(\prod_{\q \in \MM} \q\right)^{-1} = \frac{\NN(\b)}{\NN(\d_1(A))} \sum_{\MM \subseteq \QQ'}  \prod_{\q \in \MM} (-\NN(\q)^{-1}) \\
& = \frac{\NN(\b)}{\NN(\d_1(A))} \prod_{\q \in \QQ'} (1 - \NN(\q)^{-1}).
\end{align*}
Since $\tabs{M} = \mu_\a(A)$ according to Corollary \ref{elemmua} and $\b = \d_2(A) \a^{-1}$, the proof is complete.
\end{proof}

The just proved formula will be applied in the following

\begin{example}
Let $\o = \Z$ and $A = \tmatr{1}{0}{0}{4}$ as well as $\a = 2 \Z$. Theorem \ref{tmua} then yields
\[
\mu_\a(A) = \frac{\NN(4 \Z)}{\NN(2 \Z) \NN(\Z)} \prod_{\substack{\p \text{ prime ideal} \\ \p \divides 2 \Z}} (1 - \NN(\p)^{-1}) = \frac{4}{2 \cdot 1} \left(1 - \frac{1}{2}\right) = 1,
\]
which corresponds to the results of Example \ref{exz}, where we had exactly one representative of type $\tmatr{2}{*}{0}{*}$, namely $\tmatr{2}{1}{0}{2}$.
\end{example}

Since Theorem \ref{tmua} is just an intermediate result, more interesting cases than $\o = \Z$ will not be discussed at this point.

The formula for $\mu_\a(A)$ given in Theorem \ref{tmua} has several applications. Later we will see how it can be used to prove a reduction theorem in the context of Hecke algebras, but for now we will stick to the already announced goal of a formula for the number of right cosets contained in a given double coset.

\begin{theorem}\label{rcindc}
Let $A \in I$. Then
\[
\mu(A) = \NN(\f_2(A)) \prod_{\substack{\p \text{ prime ideal} \\ \p \divides \f_2(A)}} (1 + \NN(\p)^{-1}).
\]
\end{theorem}
\begin{proof}
To calculate $\mu(A)$, we have to sum over all $\mu_\a(A)$. Then we use Theorem \ref{tmua} and rewrite the obtained sum to use $\a' = \d_2(A) \d_1(A)^{-1} \a^{-1}$ as summation index:
\begin{align*}
\mu(A) & = \sum_{\a \text{ ideal in } \o} \mu_\a(A) \\
& = \sum_{\d_1(A) \divides \a \divides \d_2(A) \d_1(A)^{-1}} \frac{\NN(\d_2(A))}{\NN(\a) \NN(\d_1(A))} \prod_{\substack{\p \text{ prime ideal} \\ \p \divides \a \d_1(A)^{-1} + \d_2(A) \d_1(A)^{-1} \a^{-1}}} (1 - \NN(\p)^{-1}) \\
& = \sum_{\o \divides \a' \divides \f_2(A)} \NN(\a') \prod_{\substack{\q \text{ prime ideal} \\ \q \divides (\a')^{-1} \f_2(A) + \a'}} (1 - \NN(\q)^{-1}).
\end{align*}
Using this equality, we can prove the theorem by showing that
\[
S(\b) := \sum_{\o \divides \a \divides \b} \NN(\a) \prod_{\substack{\q \text{ prime ideal} \\ \q \divides \a^{-1} \b + \a}} (1 - \NN(\q)^{-1}) = \NN(\b) \prod_{\substack{\q \text{ prime ideal} \\ \q \divides \b}} (1 + \NN(\q)^{-1})
\]
holds for every ideal $\b$ in $\o$ (since $\b = \f_2(A)$ yields the assertion). We carry out an induction on the number of prime ideals dividing $\b$. The initial case $\b = \o$ is obvious, so we now assume that there exists a prime ideal $\p$ which divides $\b$. Write $\b = \p^m \r$ with $\p \ndivides \r$. Analogously split up every $\a$ as product of a power of $\p$ and a rest not divided by $\p$. Introducing the set $\QQ_\c$ of prime ideals dividing $\c^{-1} \b + \c$, we then have
\[
S(\b) = \sum_{\o \divides \a \divides \b} \NN(\a) \prod_{\q \in \QQ_\a} (1 - \NN(\q)^{-1}) = \sum_{k=0}^m \sum_{\o \divides \c \divides \r} \NN(\p^k \c) \prod_{\q \in \QQ_{\p^k \c}} (1 - \NN(\q)^{-1}).
\]
If $\p \ndivides \c$ and $\q$ is a prime ideal in $\o$, the definition of $\QQ_{\p^k \c}$ yields 
\[
\q \in \QQ_{\p^k \c} \quad\Leftrightarrow\quad (\q = \p \text{ and } 1 \leq k < m) \text{ or } (\q \neq \p \text{ and } \q \in \QQ_c).
\]
Using this equivalence in the above expression for $S(\b)$, by splitting up the outer sum we obtain
\begin{align*}
S(\b) = & \sum_{k=1}^{m - 1} \sum_{\o \divides \c \divides \r} \NN(\p)^k \NN(\c) (1 - \NN(\p)^{-1}) \prod_{\q \in \QQ_\c} (1 - \NN(\q)^{-1}) \\
& + \sum_{\o \divides \c \divides \r} \NN(\c) \prod_{\q \in \QQ_\c} (1 - \NN(\q)^{-1}) + \sum_{\o \divides \c \divides \r} \NN(\p)^m \NN(\c) \prod_{\q \in \QQ_\c} (1 - \NN(\q)^{-1}).
\end{align*}
Since the double sum on the right hand side is a telescoping sum, the equation simplifies to
\begin{align*}
S(\b) & = \sum_{\o \divides \c \divides \r} \NN(\p)^{m - 1} \NN(\c) \prod_{\q \in \QQ_\c} (1 - \NN(\q)^{-1}) + \sum_{\o \divides \c \divides \r} \NN(\p)^m \NN(\c) \prod_{\q \in \QQ_\c} (1 - \NN(\q)^{-1}) \\
& = \NN(\p)^m (1 + \NN(\p)^{-1}) \sum_{\o \divides \c \divides \r} \NN(\c) \prod_{\q \in \QQ_\c} (1 - \NN(\q)^{-1}) \\
& = \NN(\p)^m (1 + \NN(\p)^{-1}) S(\r).
\end{align*}
Applying the induction hypothesis, we have
\begin{align*}
S(\b) & = \NN(\p)^m (1 + \NN(\p)^{-1}) S(\r) = \NN(\p)^m (1 + \NN(\p)^{-1}) \NN(\r) \prod_{\substack{\q \text{ prime ideal} \\ \q \divides \r}} (1 + \NN(\q)^{-1}) \\
& = \NN(\b) \prod_{\substack{\q \text{ prime ideal} \\ \q \divides \b}} (1 + \NN(\q)^{-1}),
\end{align*}
which completes the proof.
\end{proof}

In the following examples Theorem \ref{rcindc} is applied in a case where $\o$ is not a principal ideal domain.

\begin{examples}\label{exmu}
Let $\o = \Z + \Z \omega$ where $\omega = \sqrt{-5}$ and $A = \tmatr{1}{0}{0}{3}$. Since $3 \o = \p_1 \p_2$ with $\p_1 = 3 \o + (\omega + 1) \o$ and $\p_2 = 3 \o + (\omega + 2) \o$, where $\p_1$ and $\p_2$ are prime ideals of norm $3$ in $\o$, Theorem \ref{rcindc} yields
\[
\mu(A) = \NN(3 \o) \prod_{\substack{\p \text{ prime ideal} \\ \p \divides 3 \o}} (1 + \NN(\p)^{-1}) = 9 \left(1 + \frac{1}{3}\right) \left(1 + \frac{1}{3}\right) = 16.
\]
Since $2 \o$ has the prime ideal decomposition $(2 \o + (\omega + 1) \o)^2$, one similarly obtains $\mu(\tmatr{1}{0}{0}{2}) = 6$. Possible choices for the six representatives are calculated in Example \ref{exdecomp}.

The above examples can be generalised: If $\o$ is a quadratic number field and $p$ a rational prime, we have $\mu(\tmatr{1}{0}{0}{p}) = (p + 1)^2$ if $p$ is split and $\mu(\tmatr{1}{0}{0}{p}) = p (p + 1)$ otherwise.
\end{examples}

To complete this section, already existing results similar to Theorem \ref{rcindc} are shortly reviewed in the following

\begin{remark}
In the case $\o = \Z$ we have a so called rationality theorem for abstract Hecke algebras with respect to unimodular groups (see e.\,g.\ \cite{Krieg:Hecke} Theorem V\,(9.3)). The proof of this theorem uses the fact that the double coset $\GL_n(\Z) P_j \GL_n(\Z)$ where $P_j$ is a diagonal matrix with $j$ diagonal entries equal to $p$ (for a fixed rational prime $p$) and the other diagonal entries equal to $1$ decomposes into exactly
\[
p^{-\frac{j (j + 1)}{2}} \sum_{1 \leq v_1 < \cdots < v_j \leq n} p^{v_1 + \cdots + v_j}
\]
right cosets with respect to $\GL_n(\Z)$. (One easily checks that for $n = 2$ and $j \in \{0, 1, 2\}$ this yields the same values for $\mu(P_j)$ as Theorem \ref{rcindc}.)

Another similar theorem does not count right cosets in double cosets but right cosets in the set of all matrices with the same determinant (modulo units). According to \cite{Newman:Integral} Theorem II.4, the set $\{A \in \o^{n \times n} \:|\: \det A \in d \o^*\}$ decomposes into exactly
\[
\prod_{\substack{p \in P \\ p \divides d}} \prod_{j=1}^{n - 1} \frac{\NN(p)^{\v_{p \o}(d) + j} - 1}{\NN(p)^j - 1}
\]
right cosets with respect to $\GL_n(\o)$ (where $P$ denotes a system of representatives of prime elements in $\o$ modulo $\o^*$).
\end{remark}

\section{Applications to congruence subgroups}\label{applcong}

In this section, an application of Theorem \ref{rcindc} to the calculation of indexes of certain congruence subgroups is presented.

\begin{corollary}\label{congindex}
Let $m \in \o$ with $m \neq 0$ and $U^0[m] = \{\tmatr{a}{b}{c}{d} \in U \:|\: b \in m \o\}$. Then the index of $U^0[m]$ in $U$ can be calculated by
\[
[U : U^0[m]] = \NN(m \o) \prod_{\substack{\p \text{ prime ideal} \\ \p \divides m}} (1 + \NN(\p)^{-1}).
\]
\end{corollary}
\begin{proof}
Let $A = \tmatr{1}{0}{0}{m}$. A simple calculation using $A \tmatr{a}{b}{c}{d} A^{-1} = \tmatr{a}{m^{-1} b}{m c}{d}$ shows that $U \cap A^{-1} U A = U^0[m]$. Since $[U : U \cap A^{-1}  U A] = \mu(A)$ (see e.\,g.\ \cite{Andrianov:Quadratic} Lemma 3.1.2), the assertion immediately follows from Theorem \ref{rcindc}.
\end{proof}

\begin{remark}
Corollary \ref{congindex} generalises a similar formula for the index of congruence subgroups in $\operatorname{SL}_2(\Z)$ which is of interest in the theory of modular forms (see e.\,g.\ \cite{Diamond:First} Section 1.2). However, such index formulae are studied also in other contexts (see e.\,g.\ \cite{Appel:Index}).
\end{remark}

\section{Applications to Hecke algebras}\label{applhecke}

As it has already been mentioned, Theorem \ref{rcindc} has been developed with the theory of Hecke algebras in mind. The applications in this field will be presented here.

Denote by $H$ the complex vector space spanned by $\{1_{U A U} \:|\: A \in I\}$ where $1_M: G \to \{0, 1\}$ is the characteristic function of the set $M$. For $A, A_1, \ldots, A_k, B, B_1, \ldots, B_m \in I$ with $U A U = U A_1 \cup \cdots \cup U A_k$ and $U B U = U B_1 \cup \cdots \cup U B_m$ where the unions are pairwise disjoint define
\[
1_{U A U} * 1_{U B U} = \sum_{k=1}^k \sum_{j=1}^m 1_{U A_i B_j}
\]
and extend this operation bilinearly to a (well-defined(!)) operation on $H$. The obtained algebra is called an (abstract) Hecke algebra; for details see e.\,g.\ \cite{Krieg:Hecke}. The formula
\[
(1_{U A U} * 1_{U B U})(C) = \tabs{\{(i, j) \:|\: A_i B_j \in U C,\; 1 \leq i \leq k,\; 1 \leq j \leq m\}},
\]
which can be found in \cite{Krieg:Hecke} I.4.4, immediately yields an algorithm for the calculation of $1_{U A U} * 1_{U B U}$.

\begin{algorithm}\label{algprod}
input: $A, B \in I$; output: $D \subseteq I$ and $c_C \in \N$ for every $C \in D$ such that
\[
1_{U A U} * 1_{U B U} = \sum_{C \in D} c_C 1_{U C U}
\]

\begin{enumerate}[(1)]
\item Decompose $U A U$ and $U B U$ into pairwise disjoint right cosets $U A_1, \ldots, U A_k$ and $U B_1, \ldots, U B_m$, respectively.
\item Let $D = \emptyset$.
\item For every pair $(i, j)$ with $1 \leq i \leq k$ and $1 \leq j \leq m$ test whether there exists a $C' \in D$ with $U A_i B_j U = U C' U$; if this is not the case, add the element $A_i B_j$ to $D$ and set $c_{A_i B_j} = 1$; otherwise, if additionally $U C' = U A_i B_j$ is fulfilled, increase $c_{C'}$ by $1$.
\end{enumerate}
\end{algorithm}

For the execution of this algorithm, a right coset decomposition of $U A U$ and $U B U$ has to be constructed explicitly in step (1). Using Theorem \ref{rcindc} we can give an algorithm that carries out this task.

\begin{algorithm}\label{algdecomp}
input: $A \in I$ and an enumeration $(Q_n)_{n \in \N}$ of $U$; output: right transversal $R$ of $U \setminus U A U$

\begin{enumerate}[(1)]
\item Calculate $k = \mu(A)$ (using Theorem \ref{rcindc}).
\item Set $R = \{A\}$ and $n = 1$.
\item\label{loop} If there exists no $B \in R$ with $U A Q_n = U B$, add the element $A Q_n$ to $R$.
\item If $\tabs{R} < k$, increase $n$ by $1$ and go back to (\ref{loop}), otherwise stop.
\end{enumerate}
\end{algorithm}

\begin{remark}\label{rrestrict}
In order to implement Algorithm \ref{algdecomp}, we have to enumerate all elements of $U$, which might not be feasible. To avoid this problem, one can use random elements instead of enumerated elements for $Q_n$. Then Algorithm \ref{algdecomp} is turned into a probabilistic algorithm which produces the desired output if it terminates. The remaining problem of the generation of random unimodular matrices will not be discussed here but is delegated to SAGE (\cite{SAGE}).
\end{remark}

Using Algorithm \ref{algdecomp} and Remark \ref{rrestrict}, we can calculate some

\begin{examples}\label{exdecomp}
Let $\o = \Z + \Z \omega$ for $\omega = \sqrt{-5}$. For $A = \tmatr{1}{0}{0}{2}$ the probabilistic decomposition algorithm terminates after an average of $14$ loop cycles and yields for example
\[
\left\{\matr{1}{0}{0}{2}, \matr{1}{1}{0}{2}, \matr{1}{\omega}{0}{2}, \matr{1}{1 + \omega}{0}{2}, \matr{2}{0}{0}{1}, \matr{2}{0}{1 + \omega}{1}\right\}
\]
as a system of representatives of $U \backslash U A U$ (with $6$ elements according to Examples \ref{exmu}). With this transversal it is then possible to use Algorithm \ref{algprod} to calculate $1_{U A U} * 1_{U A U}$; one obtains
\[
1_{U A U} * 1_{U A U} = 1_{U A_1 U} + 6 \cdot 1_{U A_2 U} + 1_{U A_3 U}
\]
with $A_1 = \tmatr{1}{0}{0}{4}$ and $A_2 = \tmatr{2}{0}{0}{2}$ as well as $A_3 = \tmatr{2}{1 + \omega}{0}{2}$.

In order to obtain a feeling for the complexity of the decomposition algorithm (a detailed analysis has to take into account the strategy for choosing the elements of $U$ and will not be carried out in this paper), we execute this algorithm for some more $A$ and obtain the following table:
\[
\renewcommand{\arraystretch}{1.3}
\begin{array}{cccrr}
\toprule
A & \d_1(A) & \d_2(A) & \mu(A) & \text{avg.\ loop cycles} \\
\midrule
\tmatr{1}{0}{0}{2} & \o & 2 \o & 6 & 14 \\
\tmatr{1}{0}{0}{1 + \omega} & \o & (\omega + 1) \o & 12 & 39 \\
\tmatr{1}{0}{0}{3} & \o & 3 \o & 16 & 53 \\
\tmatr{2}{1}{0}{2} & \o & 4 \o & 24 & 110 \\
\tmatr{\omega}{1}{0}{\omega} & \o & 5 \o & 30 & 130 \\
\tmatr{\omega}{1}{1}{2} & \o & (1 + 2 \omega) \o & 32 & 124 \\
\tmatr{\omega}{0}{0}{2} & \o & 2 \omega \o & 36 & 171 \\
\bottomrule
\end{array}
\]
\end{examples}

Theorem \ref{rcindc} can not only be used for algorithmic calculations in Hecke algebras; it has an application in the proof of a theoretical result on abstract Hecke algebras, too. In the ``classic'' Hecke algebra $H_n$ related to $\GL_n(\Z)$, certain products in $H_n$ can be reduced to products in $H_{n - 1}$ (see e.\,g.\ \cite{Krieg:Hecke} Lemma V\,(8.3)). In the case $n = 2$, the presented result yields a lemma which then leads to a generalisation of this reduction theorem to Hecke algebras related to arbitrary norm-finite Dedekind domains.

\begin{lemma}\label{prodmuo}
For all $f \in H$ define $\mu_\o(f) = \sum_{A \in R} f(A) \cdot \mu_\o(A)$ where $R$ is a system of representatives of $U \backslash I / U$. Then $\mu_\o(f * g) = \mu_\o(f) \mu_\o(g)$ for all $f, g \in H$.
\end{lemma}
\begin{proof}
By the definition of $\mu_\o$ and $*$ it suffices to prove the assertion for $f = 1_{U A U}$ and $g = 1_{U B U}$ where $A, B \in I$. For $C \in I$ with $(1_{U A U} * 1_{U B U})(C) \neq 0$ we have $C \in U A U B U$ by the definition of $*$, and \cite{Ensenbach:Determinantal} Theorem 3.1 yields $\d_1(A) \d_1(B) \divides \d_1(C)$ and thus $\d_1(A) \d_1(B) \divides \g(C)$. This implies that in the case $\d_1(A) \neq \o$ or $\d_1(B) \neq \o$ both sides of the equation $\mu_\o(1_{U A U} * 1_{U B U}) = \mu_\o(1_{U A U}) \mu_\o(1_{U B U})$ evaluate to zero, so it remains to analyse the case $\d_1(A) = \o = \d_1(B)$. In this case let $A_1, \ldots, A_k$ and $B_1, \ldots, B_m$ be systems of representatives of $U \backslash U A U$ and $U \backslash U B U$, respectively, where the $B_j$ with $\g(B_j) = \o$ have the form $\tmatr{1}{*}{0}{*}$ (without loss of generality due to Corollary \ref{unimodgen}). Let $1 \leq i \leq k$ and $1 \leq j \leq m$. If $\g(B_j) \neq \o$, then $\g(A_i B_j) \neq \o$ since the first column of $A_i B_j$ consists of linear combinations of entries of the first column of $B_j$. If $\g(B_j) = \o$, then the special structure of $B_j$ yields that the first column of $A_i B_j$ equals the first column of $A_i$. So we have $\g(A_i B_j) = \o$ if and only if $\g(A_i) = \o$ and $\g(B_j) = \o$. Since according to the definition of $\mu_\o$ and $*$ we have
\[
\mu_\o(1_{U A U} * 1_{U B U}) = \tabs{\{(i, j) \:|\: \g(A_i B_j) = \o,\; 1 \leq i \leq k,\; 1 \leq j \leq m\}},
\]
the just proved characterisation of $\g(A_i B_j) = \o$ used to split up the right hand side as product of two cardinalities yields the assertion.
\end{proof}

Now the desired reduction theorem can be stated and proved.

\begin{theorem}\label{reduction}
Let $a, b, c \in \o$ and $A = \tmatr{1}{0}{0}{a}$, $B = \tmatr{1}{0}{0}{b}$ as well as $C = \tmatr{1}{0}{0}{c}$. Then $(1_{U A U} * 1_{U B U})(C) = 1$, if $c \in a b \o^*$, and $(1_{U A U} * 1_{U B U})(C) = 0$ otherwise.
\end{theorem}
\begin{proof}
With $R$ as in Lemma \ref{prodmuo} write
\begin{align*}
\mu_\o(1_{U A U} * 1_{U B U}) & = \sum_{D \in R} (1_{U A U} * 1_{U B U})(D) \cdot \mu_\o(D) \\
& = \sum_{\substack{D \in R \\ D \notin U A B U}} (1_{U A U} * 1_{U B U})(D) \cdot \mu_\o(D) + (1_{U A U} * 1_{U B U})(A B) \cdot \mu_\o(A B).
\end{align*}
Using Lemma \ref{prodmuo}, Theorem \ref{tmua} and the multiplicity of the norm, we have
\[
\mu_\o(1_{U A U} * 1_{U B U}) = \mu_\o(1_{U A U}) \mu_\o(1_{U B U}) = \NN(\d_2(A)) \NN(\d_2(B)) = \NN(\d_2(A B)) = \mu_\o(A B),
\]
so
\[
\sum_{\substack{D \in R \\ D \notin U A B U}} (1_{U A U} * 1_{U B U})(D) \cdot \mu_\o(D) + (1_{U A U} * 1_{U B U})(A B) \cdot \mu_\o(A B) = \mu_\o(A B).
\]
Since all numbers in this equation are non-negative integers and $(1_{U A U} * 1_{U B U})(A B) \geq 1$ by the definition of $*$, we have $(1_{U A U} * 1_{U B U})(A B) = 1$ and $(1_{U A U} * 1_{U B U})(D) \cdot \mu_\o(D) = 0$ for all $D \in R$ with $D \notin U A B U$. Since $\mu_\o(C) \geq 1$ as $\g(C) = \o$, these equations imply $(1_{U A U} * 1_{U B U})(C) = 0$ if $C \notin U A B U$ and $(1_{U A U} * 1_{U B U})(C) = 1$ if $C \in U A B U$, where the latter condition is equivalent to $c \in a b \o^*$, which proves the assertion.
\end{proof}

\bibliographystyle{unsrt}
\bibliography{countcosets}

\end{document}